# On Path Decomposition Conjecture of Tibor Gallai

Dhananjay P. Mehendale
Sir Parashurambhau College, Tilak Road, Pune 411030,
India

## Abstract

We settle the Path Decomposition Conjecture (P.D.C.) due to Tibor Gallai for minimally connected graphs, i.e. trees. We use this validity for trees and settle the P. D. C. using induction on the number of edges for all connected graphs. We then obtain a new bound for the number of paths in a path cover in terms of the number of edges using idea of associating a tree with a connected graph. We then make use of a spanning tree in the given connected graph and its associated basic path cover to settle the conjecture of Tibor Gallai in an alternative way. Finally, we show the existence of Hamiltonian path cover satisfying Gallai bound for complete graphs of even order and discuss some of its possible ramifications.

1. **Introduction:** All graphs that we consider here are simple, finite, and undirected. Also, they are connected graphs unless mentioned otherwise. By [x] we denote the integral part of x. If a connected graph G has a vertex, v, with degree, $d$ say, then (it is easy to check that) at least $\left[\frac{(d+1)}{2}\right]$ paths must pass through v. Note further that if the degree of vertex v is even then vertex v may or may not be end vertex for any path among the $\left[\frac{(d+1)}{2}\right]$ paths that pass through it, but if the degree of vertex v is odd then there must be at least one path among the paths passing through v for which v is end vertex for that path. (e.g. for a vertex of degree 3 there must be at least two paths to which these edges belong and between these paths the vertex under consideration should be end vertex for at least one path.)

    Let $d_{max}$ be the maximum degree for vertices of G then in the disjoint path cover for G there should be at least $\left[\frac{(d_{max}+1)}{2}\right]$ paths

in it, i.e. $\left\lceil \dfrac{(d_{max}+1)}{2} \right\rceil$ forms a lower bound on the number of paths in the path cover of G. What is the upper bound on the number of paths in its disjoint path cover? For a connected graph G containing n vertices the conjecture due to T. Gallai (1962), popularly known as Path Decomposition Conjecture (P. D. C.), asserts that G can be decomposed into at most $\left\lceil \dfrac{(n+1)}{2} \right\rceil$ disjoint paths. In the same year, it has been shown by L. Lovasz [1] that P. D. C. is true for graphs in which all degrees are odd. Further, let H be a subgraph of G induced by vertices of even degree in G, then it has been shown by L. Lovasz that the conjecture is true if H contains at most one vertex. Pyber [2] extended Lovasz's result and proved that the conjecture is true if H is a forest. These results were further extended by G. Fan [3]. Fan proved that the conjecture is true if H can be obtained from the empty set by a series of so-defined "alpha" operations. As a corollary, Fan showed that the conjecture is true if each block of H is a triangle-free graph of maximum degree at most 3.

      Our approach to this problem in this paper is **new and constructive**. We first establish the validity of Gallai's conjecture for trees, i.e. the minimally connected graphs and proceed by induction on the number of edges to settle Gallai's conjecture for all connected graphs. We then proceed to obtain a new upper bound for the number of paths required in the path decomposition of a connected (n, e) graph G in terms of its edge number, e. We then proceed with a new algorithmic approach and establish Gallai's conjecture using idea of B. P. D. for spanning tree of a connected graph and extension of its basic paths into trees containing maximally extended paths that form a cover for connected graph under consideration satisfying Gallai's bound. Finally, we proceed to prove Gallai's conjecture for even ordered maximally connected graphs, i.e. complete graphs of even order. We show something more in this case, namely, the existence of Hamiltonian path decomposition, i.e. decomposition in terms of Hamiltonian paths, for complete graphs of even order and indicate some of the possible outcomes of this approach.

2. **Gallai's Conjecture for Trees:** We begin with proving the conjecture for minimally connected graphs, i.e. trees. We give three proofs for this case. The first proof is direct and elementary while the second proof makes use of the result of Pyber [2] mentioned above, and the third very elegant proof is due to Riko Winterle.

**Definition 2.1:** A path P in a connected graph G is called **maximally extended path** if its extension at either end by appending an edge is not possible, i.e. if such an extension by an edge leads to formation of a circuit or subcircuit in G. In other words, a path P in a connected graph G is called maximally extended path if there is no new vertex other than the one that already belongs to P which is adjacent to one of the two end vertices of P.

**Definition 2.2:** A Path P connecting two pendant vertices in a tree T is called **peripheral path** if the removal of P (i.e. the removal of edges of P) from tree T leads to another tree T*, i.e. T - P = T*.

**Definition 2.3:** A Path P connecting two vertices in a graph G is called **peripheral path** if the removal of P from tree G leads to another connected graph G*, i.e. if G - P = G*, then G* is connected.

**Definition 2.4:** A Path P connecting two vertices in a graph G is called **maximally extended peripheral path (m.e.p.p.)** if it is both, maximally extended as well as peripheral.

**Remark 2.1:** A peripheral path P having two pendant vertices in a tree T as its end vertices is maximally extended.

**Remark 2.2:** There can exists many distinct disjoint path covers, each made up of a set of certain "maximally extended peripheral paths (m.e.p.p)" made up of paths of different lengths for a given connected graph. As there is no restriction on the staring vertex or edge to begin construction of m.e.p.p.(s) in succession, therefore, only a part of some m.e.p.p. in certain path cover may appear as m.e.p.p. in some other path cover. It depends on where we start and how we proceed with the formation of certain m.e.p.p. in the construction procedure of such paths in succession. The appearance of only a part of some m.e.p.p. in certain path cover as m.e.p.p. in some other path cover occurs because the edges which were present in the earlier m.e.p.p. in the first cover are nonexistent as they have become part of some other m.e.p.p. formed earlier in the second path cover.

**Lemma 2.1:** Every tree contains a maximally extended peripheral path.

**Proof:** We proceed by induction on n, the number of vertices in tree T.

   1) It is easy to check the result for n = 1,2,3,4 etc.

   2) Now, let T be a tree with n >4 vertices. Remove a leaf in tree T and suppose it leads to tree T1. By induction T1 contains a maximally extended peripheral path P. Now attach back the removed leaf. It is clear to see that P still remains maximally extended peripheral in T.

                                                               □

**Lemma 2.2:** When a maximally extended peripheral path is removed from a tree T leading to a tree T* then the count of number of vertices in T* is ≤ (n-2), i.e. the count of the number of vertices in T* is less than the count of vertices in T by at least two.

Proof: When a maximally extended peripheral path P is removed at least the pendent vertices contained in it get removed assuredly.

                                                               □

**Theorem 2.1** Every tree satisfies Gallai's conjecture.

**First Proof:** We proceed by induction on n, the number of vertices in tree, T, say.

   1) It is easy to check the result for n = 1,2,3,4 etc.

   2) Now, let T be a tree with n (>4) vertices. Remove a maximally extended peripheral path P in tree T and suppose it leads to tree T*. By induction T* satisfies Gallai's conjecture, therefore, the path decomposition of T* requires $\left\lceil \frac{(n-1)}{2} \right\rceil$ paths. Now attach back the removed maximally extended peripheral path P. It is clear that T therefore requires $\left\lceil \frac{(n-1)}{2} \right\rceil + 1$ paths in its path cover.

Now, since $\left\lceil \frac{(n-1)}{2} \right\rceil + 1 \leq \left\lceil \frac{(n+1)}{2} \right\rceil$, the result holds for tree T.

                                                               □

**Second Proof:** It is clear to see that every subgraph of tree is either a tree or a forest. Therefore, a subgraph of a tree T induced by vertices of even degree in T will be essentially a forest, and so, as per the result of Pyber [2] the conjecture due to T. Gallai will hold good for a tree T.

□

**Third Proof (Riko Winterle):** Let G be a connected (n, e) graph and let C be a minimal covering for G by edge disjoint "paths of special kind" having strictly different starting and ending vertices but these paths may contain a cycle if they intersect themselves. Now, these paths cannot have common end vertices because otherwise we can join those paths with same end vertices to get a covering D which will contain fewer number of such paths, a contradiction, since C is assumed to be minimal. Hence if cardinality of paths in C is P then the cardinality of distinct end vertices is 2p. Now, since number of vertices in G is equal to n, we have $2p \leq n$, therefore, $p \leq \left\lceil \frac{(n+1)}{2} \right\rceil$. Now, if graph G is some tree T then the "paths of special kind" mentioned above will be just the ordinary paths!

□

**Conjecture 2.1(An equivalent of Gallai's conjecture):** Let G be a connected (n, e) graph. Let $l$ denote the average length of maximally extended peripheral paths obtained in a sequence that covers G. Gallai's conjecture asserts that $l \geq \left\lceil \frac{(e+1)}{\left\lceil \frac{(n+1)}{2} \right\rceil} \right\rceil$. In other words, if $l < \left\lceil \frac{(e+1)}{\left\lceil \frac{(n+1)}{2} \right\rceil} \right\rceil$ then G must be disconnected.

We now proceed to settle the Gallai's conjecture for connected graphs. We proceed by induction on number of edges and show that the result holds in general.

**Lemma 2.3:** Let G be a connected (n, e) graph, i.e. G contains n vertices and e edges. Let u, v be some vertices of G and edge uv does not belong to G. Then one can build a path decomposition for G in terms of m.e.p.p.s, C say, such that if the total number of paths in this path decomposition of G is equal to $\left[\frac{(n+1)}{2}\right]$ then there exists an m.e.p.p. in C with vertex u (or v) as end vertex and this m.e.p.p. in C does not contain vertex v (or u).

**Proof:** Suppose not. It then follows that for every path decomposition we have paths for which either
(i) Both u and v appear as intermediate vertices (i.e. vertices of degree two) but this will contradict the total possible degree for both the vertices u and v in a graph containing n vertices.
(ii) One vertex u as end vertex and other vertex v as intermediate vertex (i.e. vertex u with degree one and vertex v with degree two) but this will contradict the total possible degree for vertex v in a graph containing n vertices.
(iii) Both u and v appear as end vertices in every path. But in this case if all the paths are of full length that is possible for a graph on n vertices then it will contradict the total possible degree for some vertex other than u and v in a graph containing n vertices. So, in such case the lengths of such paths must go on reducing and as stated in remark 2.2 we can reshuffle the m.e.p.p.s to bring vertex v in some path to the intermediate position by extending a path terminated at v along some other path starting u and ending at v and create a path starting at u and not containing v.
□

**Theorem 2.2(Gallai's Conjecture) :** There exists a path decomposition for every connected graph on n vertices and e edges made up of maximally extended peripheral paths, m.e.p.p.(s), containing in all $\leq \left[\frac{(n+1)}{2}\right]$ disjoint paths.

**Proof:** As proved in above the result holds for minimally connected graphs on n vertices, i.e. for connected graphs on n vertices and (n-1) edges, thus the first step of induction is clear.
We assume by induction that there exists desired path decomposition formed by m.e.p.p.(s) for all connected graphs containing n vertices and (e-1) edges and proceed to show that it holds for all graphs containing n vertices and e edges.

Let G be a connected (n, e) graph and let edge uv belongs to edge set of G. Let G' = G – uv, i.e. G' is an edge deleted subgraph of G obtained by deleting edge uv from it. We assume by induction the result for G' and show that it holds for G.

Now, if the number of m.e.p.p.s required in the path decomposition are strictly less than $\left\lceil \frac{(n+1)}{2} \right\rceil$ then we can take uv as an additional path and we have a path decomposition containing m.e.p.p.s $\leq \left\lceil \frac{(n+1)}{2} \right\rceil$. If on the other hand the number of m.e.p.p.s required in the path decomposition are exactly equal to $\left\lceil \frac{(n+1)}{2} \right\rceil$ then the lemma 2.3 assures us of the existence of path decomposition for G' in terms of m.e.p.p.s, and further it assures the existence of an m.e.p.p. which has u (or v) as end vertex and which does not contain v (or u). So, we can append edge uv to this m.e.p.p. and again have a path decomposition for G in terms of m.e.p.p.s $\leq \left\lceil \frac{(n+1)}{2} \right\rceil$. Hence etc.

∎

3. **Associating a Tree with a Connected Graph:** In this section we introduce an idea associating a tree with a connected graph which can be useful at some other place as well. We show here its use to obtain a new upper bound for the number of paths required in the path decomposition of a connected (n, e) graph G in terms of its edge number, e. We associate with (or, look at a) connected graph (as) a rooted, unordered, pseudo tree. To describe this idea let us begin with an example:

**Example 3.1:** Consider following connected graph, G say, and its associated rooted, unordered pseudo tree, T say:

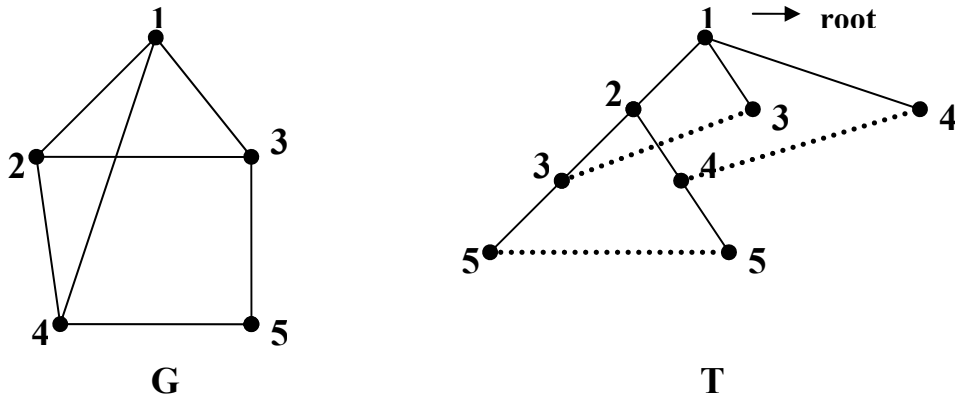

G          T

This tree is called rooted because it has a root (vertex with label 1), unordered because it contains multiple vertices with same label, and pseudo because it contains pseudo edges joining vertices with same label shown by dotted line segments which are actually nonexistent.

It is easy to construct a conveniently labeled (this is actually not essential but convenient) graph from given unlabeled graph and further its associated rooted, unordered, pseudo tree with the help of the following two procedures:

**Graph Labeling Procedure:**

A) Choose any vertex in the given unlabeled graph to label it as "1".
B) Assign labels "2", "3", …., "d1+1" to the vertices adjacent to vertex with label "1".
C) Take vertex with label "2" and assign labels "d1+1", "d1+2", ….to vertices adjacent to vertex with label "2" other than already labeled adjacent vertices.
D) Continue taking vertices with labels "3", "4", …. and repeat the above procedure till every vertex gets a label.

**Tree Construction Procedure:**

a) Take vertex with label "1" in the above labeled graph as root (1st level)
b) Create vertices "2", "3",…. and edges "12", "13", …. (2nd level)
c) Take vertex with label "2" in the 2nd level and create new level of vertices e.g. "3", "4", ….etc. which are adjacent to vertex with label "2".
d) Continue till we reach the last level containing vertices with highest label "n" when the graph contains n vertices.
e) Join vertices with same label by a pseudo edge (dotted line segment).

**Remark 3.1:** To avoid formation of any cycle (which should not exist in a tree) we take such vertices as additional vertices in the next level so clearly the rooted tree contains more vertices than the one contained in the original graph.

**Remark 3.2:** It is straightforward to check that if the connected graph G under consideration is an (n, e) graph, i.e. it contains n vertices and e edges, then its associated rooted, unordered, pseudo tree contains (e+1) vertices.

**Remark 3.3:** It is easy to see further that in all there are (e+1-n) repetitions of vertices, i.e. there are in all (e+1-n) vertices which show appearance more than once in the associated rooted, unordered, pseudo tree. Thus, there are in all (e+1-n) number of **distinct** pairs of vertices which are identical and so each pair of identical vertices is connected by a pseudo edge (shown by a dotted line segment in the above figure).

**Theorem 3.1:** Every connected (n, e) graph can be decomposed into disjoint paths $\leq \left\lceil \frac{(e+2)}{2} \right\rceil$.

**Proof:** As seen above with every connected (n, e) graph G we can associate a (rooted, unordered, pseudo) tree containing edges of G and some pseudo edges, and this tree contain (e+1) vertices. By theorem 2.1, it is clear that this tree (ignoring the pseudo edges) can be covered by disjoint paths $\leq \left\lceil \frac{(e+2)}{2} \right\rceil$ in number. Hence etc.

□

**Algorithm 3.1:**
1) Carry out "Graph Labeling Procedure" and its associated "Tree Construction Procedure" for given (n, e) graph G and obtain the associated tree T.
2) Choose some pseudo edge, among some distinct (e+1-n) pseudo edges (which will be taken up later in a sequence for the same procedure that we are going to develop for the chosen pseudo edge), joining identical vertex pair, (k, k) say.
3) Start from label k at one end of pseudo edge proceed to select edges of G in a sequence, one adjacent to other, and build path in T till either, (i) we reach a vertex which is end vertex of some other pseudo edge among the distinct (e+1-n) pseudo edges, or (ii) the path cannot be extended (i.e. leads to formation of some subcircuit, or next adjacent edge is already used up in already formed other path).
4) Start now from label k at the other end of pseudo edge and repeat procedure in 3). We have thus completed formation of $1^{st}$ path pair.
5) Complete the procedure of formation of all distinct path pairs, (e+1-n) in number.
6) Now, form single paths, minimum possible in number, for any left out edges of G so that every edge of G now either belongs to a path in a path pair or belongs to some single path, among the paths built from left out edges of G in T after completing the procedure of formation of path pairs mentioned above.

□

**Remark 3.5:** We have thus completed the decomposition of T which will lead to decomposition of G when pseudo edges connecting path pairs will be contracted and other hanging pseudo edges are removed.

**Example 3.2:**

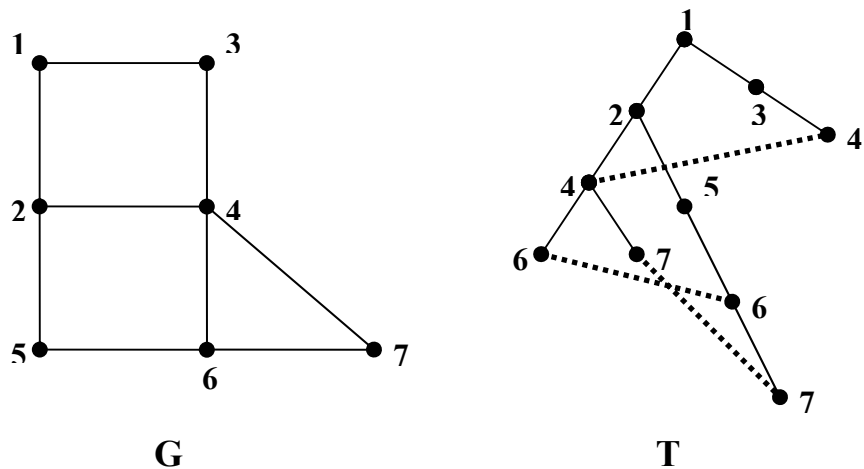

**G**            **T**

By proceeding as per above give algorithm we get following path pairs for above example:
{(74),(76)}, {(64), (65213)}, {(42), (43)}. Further, by contracting pseudo edge, these path pairs can be combined into single paths and thus we get the path decomposition for G in terms of three paths as follows: {(476), (465213), (243)}.

        We now state below a conjecture equivalent to Gallai's conjecture as follows:

**Conjecture 3.2:** For a given connected (n, e) graph G let T be its associated rooted (having root at vertex with label 1), unordered, pseudo tree. Let S be the cardinality of single paths and D be the cardinality of disjoint path pairs with identical common end vertex connected by a pseudo edge, obtained as a result of performing Algorithm 3.1 on associated T of G, then there exists at least one result of this algorithm for which

$$(D + S) \leq \left\lceil \frac{(n+1)}{2} \right\rceil.$$

4. **Spanning Tree, its Path Decomposition, and Gallai's Conjecture:** In this section we proceed to settle the Gallai's conjecture for connected graphs. We have already proved Gallai's conjecture for minimally connected graphs on n vertices, i.e. trees, in Section 2. We proceed here with an existential proof for the conjecture which also offers a constructive procedure to find the actual disjoint path cover for any connected graph.

Consider a connected (n, e) graph G. Every connected graph contains a spanning tree. So, let T be a spanning tree of G. It is clear that this tree T is (n, n-1) graph. By theorem 2.1 T can be decomposed into disjoint paths $\leq \left\lceil \frac{(n+1)}{2} \right\rceil$ in number. We by choice select the path cover for T which actually contains $\left\lceil \frac{(n+1)}{2} \right\rceil$ disjoint paths. It is easy to check that this path cover will contain (i) (k-1) number of 2-paths (paths of length two) and one 1-path (path of length one) when n = 2k, and (ii) (k-2) number of 2-paths (paths of length two) and two 1-paths (paths of length one) when n = (2k-1). We call this path decomposition **Basic Path Decomposition (B. P. D.)**. We extend these basic paths forming the basic path cover by other edges of G so that every edge of G will now belong to some tree. We show that it is possible to extend the basic paths into trees so that we can incorporate every edge of G in at least one or more trees through the extension procedure of basic paths described below. We essentially show that this procedure of starting with a basic path cover and go on growing these basic paths by other edges of G is possible to continue without increasing the count of paths, i.e. without requiring the addition of any altogether new path which is not obtained as an extension of some basic path till the consumption of all the edges of G.

In nutshell, our procedure is as follows:

**Algorithm 4.1:**

1) We take a labeled copy of G, the given connected graph.
2) We select any spanning tree T of G and take its copy.
3) We find B. P. D. for T containing $\left\lceil \frac{(n+1)}{2} \right\rceil$ disjoint paths and take a copy of this cover for T.
4) We select, one by one, paths in B. P. D. and attach all possible edges, i.e. chords, of G not in T such that no cycle gets formed and the result of this appending is a **tree**.
5) We carry out the procedure in 4) on each basic path and construct a tree from each basic path. (Note that those edges are not added in a basic path which leads to formation of a subcircuit and each addition leads strictly to formation of a tree.)
6) We associate a number called weight with each edge of G which is defined as the number of occurrences of that edge in different trees

formed from the basic paths and further construct weighted adjacency matrix associated with G.
7) From each tree among the trees formed from basic paths we form a longest possible path by selecting some edges, and we call these selected edges **used edges.** And all other edges on that tree we call **unused edges.**
8) We form these longest possible paths in such a way that we obtain now a modified weighted adjacency matrix from the weighted matrix associated with G with which we started and in which we make the following changes: 1) We make weight of every **used** edge equal to zero. Since we form only one path from one tree so, further 2) we subtract the count corresponding to the other edges on this tree which are **not used** on this path from the weight of the unused edges such that after this reduction in weight the weight of unused edges in the resulting modified weighted adjacency matrix still remains **positive**.
9) We continue the procedure in 8) till we form the longest possible paths by this procedure from every tree so that every edge of G now belongs to the category of used edges.

□

**Theorem 4.1:** Every connected graph satisfies Gallai's conjecture and algorithm 4.1 offers the desired disjoint path cover.

**Proof:** Let G be a connected (n, e) graph. As per algorithm 4.1 we choose some spanning tree T of G. We find B. P. D. of tree T made up of $\left\lceil \frac{(n+1)}{2} \right\rceil$ disjoint paths; 2-paths and 1-paths. Thus, we have **restricted** by choice the count of basic paths such that it automatically satisfies bound by the conjecture of Gallai. It is clear to see that every edge of G not belonging to any basic path can be appended to some basic path or to extended basic path or to the tree obtained by appending edges of G to basic paths. Note that each of chords, i.e. edges of G not belonging to tree T, when appended to tree T, it creates only one cycle. Therefore, there can exist only one basic path in the associated B. P. D. to which when some chord under consideration is appended will create a cycle and in case of other basic or extended paths this chord will get added as a part of some extended path or tree. As per algorithm 4.1 we extend separately these basic paths into trees, each basic path extends into a separate tree, and pick paths on these trees in a sequence made up of edges with minimum weight such that a) only the weight of **used edges** on the paths are made zero everywhere, and b) weight of **unused edges** always remains positive till it has unused status implying their presence on other trees associated with other basic paths. When all paths, one from each tree, are selected then these paths together finally contain each edge of G only once in some path since its weight is made zero amounting to its removal when it attains used status.

□

**Example 4.1:**

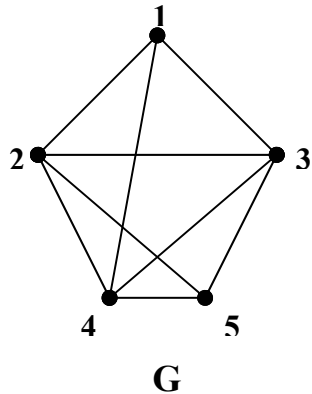

G

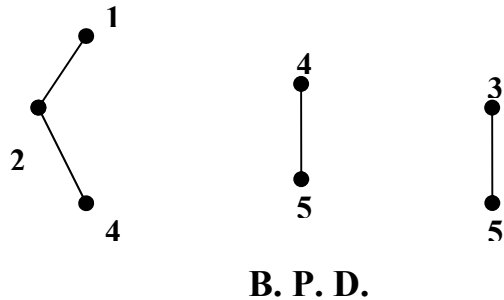

**B. P. D.**

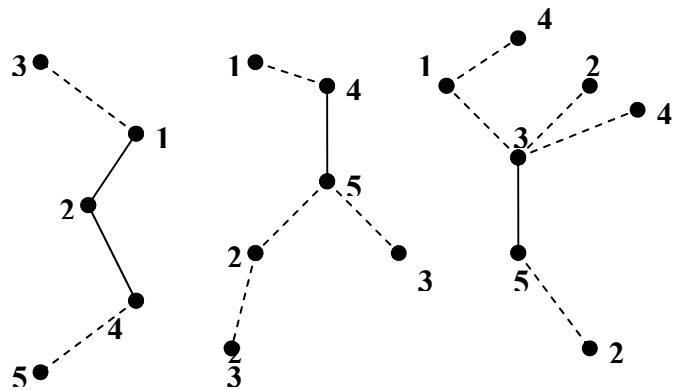

**Trees associated with Paths in B.P.D.**

From the associated trees in this example we form following paths:
$$3124$$
$$14523$$
$$435$$
as per the algorithm 4.1 given above.

5. **Hamiltonian Path Decomposition and Path Tables:** In the section 2 we proved Gallai's conjecture for minimally connected graphs, i.e. trees. In this section we proceed to prove Gallai's conjecture for even ordered maximally connected graphs, i.e. complete graphs of even order. We show something more in this case, namely, the existence of **Hamiltonian path decomposition**, i.e. decomposition in terms of Hamiltonian paths, for complete graphs of even order. We show the existence of path decomposition for $K_{2n}$ in terms of $n$ number of disjoint Hamiltonian paths. This existence of path cover in terms of Hamiltonian paths easily follows from Lovasz's theorem. We use this Hamiltonian path cover for complete graphs to construct a path cover called pseudo Hamiltonian path cover for any connected graph. We see finally that it opens up a way to construct path decomposition, for any given connected graph containing even/odd number of vertices. For the sake of easy referencing some preliminaries are in order:

A graph is called **even (odd) ordered** or simply **even (odd)** if it contains even (odd) number of vertices. We choose natural numbers as labels for its vertices. Thus, for any graph G on $n$ vertices its vertex set is:
$$V(G) = \{1, 2, ...., j, ...,n\}.$$

We denote the edges of G connecting vertex labels u and v in V(G) by (u, v), or simply by uv. Since all graphs are undirected uv and vu denote the same edge. A **path** is a sequence of distinct vertex labels which are adjacent in succession. A path starting at a vertex $u_1$ passing through intermediate vertices, $u_2, u_3, \cdots, u_{(n-1)}$ in a sequence and terminating at a vertex $u_n$ is denoted by $u_1 u_2 u_3 \cdots u_{(n-1)} u_n$. If $u_1 u_2 u_3 \cdots u_{(n-1)} u_n$ is a path then $u_1$ is adjacent to $u_2$, $u_2$ is adjacent to $u_3$,......, $u_{(n-1)}$ is adjacent to $u_n$ in G and all $u_i$, $i = 1, 2, 3, ....,n$ are distinct vertices in V(G). The vertices occurring at the two ends of path are called **end vertices.** The edges occurring at the two ends of path are called **end edges**. Those vertices that appear at intermediate position are called **intermediate vertices.** Thus, the above given path $u_1 u_2 u_3 \cdots u_{(n-1)} u_n$ contains $u_2, u_3, \cdots, u_{(n-1)}$ as intermediate vertices and vertices $u_1$ and $u_n$ as end vertices. Further, it consists of edges, $u_1 u_2, u_2 u_3, u_3 u_4, \ldots$, etc. with $u_1 u_2$ and $u_{(n-1)} u_n$ as end edges. Those edges that appear at intermediate position, so that, when they are removed the path under consideration may split into two or more paths of shorter length, are called **intermediate edges.** A **broken path** or **pseudo path** is a path which contains some edges, called **pseudo edges**, expressed by putting a cross symbol, ×, in between, which are actually absent in the path, and so, such path, though expressed as a single path is actually composed of more than one subpaths. Let path $u_1 u_2 u_3 \cdots u_{(n-1)} u_n$ be a pseudo path and let an end edge, say, $u_1 u_2$ be absent in it then we express this absentee by putting a cross × between vertex labels $u_1$ and $u_2$ to indicate that this join is artificial and there is no edge present in actuality in the graph under consideration. Thus, such an edge is expressed as $u_1 \times u_2$, and the corresponding pseudo path under consideration is expressed as $u_1 \times u_2 u_3 \cdots u_{(n-1)} u_n$. It is important to note here that such a pseudo path can still be expressed as a single path by just deleting the end vertex in this path, like, $u_2 u_3 \cdots u_{(n-1)} u_n$. Note that this has become possible because of the occurrence of broken or pseudo edge at the end of the pseudo path $u_1 \times u_2 u_3 \cdots u_{(n-1)} u_n$. Instead, if the absent or broken or pseudo edge exists as an intermediate edge, e.g. let the edge joining intermediate vertex labels $u_j$ and $u_{(j+1)}$ be absent, where $j \neq 1, (n-1)$, then such pseudo path is expressed as $u_1 u_2 u_3 \cdots u_j \times u_{(j+1)} \cdots u_{(n-1)} u_n$. Note that this pseudo path is actually made up of two paths, $u_1 u_2 u_3 \cdots u_j$ and $u_{(j+1)} \cdots u_{(n-1)} u_n$, and **as it is** it cannot be treated as a single path.

**Proposition 5.1:** There exists Hamiltonian path decomposition, consisting of n disjoint Hamiltonian paths, for $K_{2n}$.

**Proof:** It has been shown by L. Lovasz that Path decomposition conjecture due to T. Gallai is true for graphs in which all vertices have odd degrees. Now, for $K_{2n}$, the degree of every vertex is (2n-1) and so odd. Therefore, there exists path decomposition for this graph as per the Lovasz theorem containing less than or equal to $\left\lceil \frac{(2n+1)}{2} \right\rceil = n$ paths.

Now:
1) The longest length of a path in a graph on 2n vertices can be (2n-1).
2) The total number of edges in $K_{2n}$ is (n)(2n-1).
3) The maximum number of paths forming a disjoint path cover can be equal to $\left\lceil \frac{(2n+1)}{2} \right\rceil = n$.

Thus, to fulfill the constraint on length of longest path, the constraint on incorporation of every edge among the total (n)(2n-1) edges in some one path to form a decomposition, and the constraint on the maximum count of such paths due to result of Lovasz, we should have exactly n paths in the decomposition and each such path must be Hamiltonian.

□

A **Hamiltonian path table** *P* is a collection of disjoint Hamiltonian paths collected in a table which together form the disjoint path cover for the complete graph of even order.

$P$, path table for $K_{2n}$, is actually the collection of *n* number of disjoint Hamiltonian paths, as shown below:

$$P = \begin{bmatrix} 1 & 2 & \cdots & j & \cdots & 2n \\ 3 & a_2^3 & \cdots & a_j^3 & \cdots & (2n-2) \\ \vdots & \vdots & \vdots & \vdots & \vdots & \vdots \\ j & a_2^j & \cdots & a_j^j & \cdots & (2n-j+1) \\ \vdots & \vdots & \vdots & \vdots & \vdots & \vdots \\ (2n-1) & a_2^{2n-1} & \cdots & a_j^{2n-1} & \cdots & 2 \end{bmatrix}$$

**Examples:** We state Hamiltonian path tables for some complete graphs of even order:

1. For 2n = 2, $P = [1\ 2]$
2. For 2n = 4, $P = \begin{bmatrix} 1 & 2 & 3 & 4 \\ 3 & 1 & 4 & 2 \end{bmatrix}$
3. For 2n = 6, $P = \begin{bmatrix} 1 & 2 & 3 & 4 & 5 & 6 \\ 3 & 1 & 5 & 2 & 6 & 4 \\ 5 & 3 & 6 & 1 & 4 & 2 \end{bmatrix}$
4. For 2n = 8, $P = \begin{bmatrix} 1 & 2 & 3 & 4 & 5 & 6 & 7 & 8 \\ 3 & 1 & 4 & 2 & 7 & 5 & 8 & 6 \\ 5 & 2 & 6 & 1 & 7 & 3 & 8 & 4 \\ 7 & 4 & 6 & 3 & 5 & 1 & 8 & 2 \end{bmatrix}$

Note that the path table, $P$, for $K_{2n}$ is not unique and we can have its many distinct avatars! It follows from the following proposition. It should be noted that the action of permutations on the path tables considered below only the paths change and the underlying labeled graph remains intact, i.e. by this action we are essentially choosing different edge sequence in the same labeled graph as our new paths in place of old ones.

**Proposition 5.2:** A Hamiltonian path table, $P$, for $K_{2n}$ survives under the action of permutation of order 2n. In other words, let $\sigma \in S_{2n}$, where, $S_{2n}$ is group of permutations on 2n symbols, then $\sigma(P) = P^*$, where $P^*$ is again a Hamiltonian path table.

**Proof:** A permutation is one-one onto map. All the vertex pairs are present in any Hamiltonian path table, and exist in some one and only one path. Under the action of $\sigma \in S_{2n}$ every vertex pair that is present in $P$ will also remain present in $P^*$ but will only change its position. If under $\sigma \in S_{2n}$ $k \to i$ and $l \to j$ then vertex pair $(k, l) \to (i, j)$, i.e. in the place of vertex pair $(k, l)$ in $P$ there will appear vertex pair $(i, j)$ in $P^*$. Now, where will be then the vertex pair $(k, l)$? Obviously, there will be $m \to k$ and $n \to l$ defined by the $\sigma \in S_{2n}$, so the pair $(k, l)$ in path table $P$ will appear in the path table $P^*$ at the place of pair $(m, n)$, and so on.

□

A **Pseudo Hamiltonian path table** $P(G)$ associated with a connected (or disconnected) graph on $2n$ vertices is a collection of disjoint pseudo Hamiltonian paths collected in a table which is obtained by introducing cross symbols " × " between the vertex labels which are **not** adjacent in the given connected (or disconnected) graph G.

**Examples:** We state pseudo Hamiltonian path tables for some graphs of even order:

1. For 2n = 4, $P(G) = \begin{bmatrix} 1 & 2 & \times & 3 & 4 \\ 3 & 1 & \times & 4 & 2 \end{bmatrix}$

2. For 2n = 6, $P(G) = \begin{bmatrix} 1 & 2 & \times & 3 & \times & 4 & \times & 5 & 6 \\ 3 & \times & 1 & 5 & \times & 2 & 6 & \times & 4 \\ 5 & 3 & \times & 6 & 1 & \times & 4 & \times & 2 \end{bmatrix}$

Note that even if the order of graph is odd, i.e. suppose it contains ($2n$-1) vertices, still we can start with path table for $K_{2n}$ and introduce appropriate crosses representing disconnection or pseudo connection between the vertices and obtain thus pseudo Hamiltonian path table for $K_{(2n-1)}$ under consideration. For example, consider the following pseudo path table:

Let 2n = 6, and $P(G) = \begin{bmatrix} 1 & 2 & 3 & 4 & 5 & \times & 6 \\ 3 & 1 & 5 & 2 & \times & 6 & \times & 4 \\ 5 & 3 & \times & 6 & \times & 1 & 4 & 2 \end{bmatrix}$

It represents pseudo path cover for $K_6$ and the edges together represent actually the complete graph $K_5$, since we have broken only those all connections (adjacencies) which represent adjacencies of the vertex labeled as 6 in this path cover. The paths in this path cover are as given below:
$$1\ 2\ 3\ 4\ 5$$
$$3\ 1\ 5\ 2$$

5 3
1 4 2

All these paths together form a path decomposition for $K_5$. Further, this collection is not minimal as we can easily transform it to a collection of 3 paths given below by rearrangement as follows: We keep first path as it is. We take piece 31 of second path and join it to last path 142 and get new path 3142. We take 3142 and extend it further by joining it to 25 in second path to arrive at path 31425. We then pick left out part 15 of second path and join it to third path 53 to get in effect new path 153. Thus, we have new path decomposition

1 2 3 4 5
3 1 4 2 5
1 5 3

As already mentioned, in the pseudo Hamiltonian path table the paths in the rows containing at least one cross symbol representing disconnection are obviously not Hamiltonian paths but are actually made up of two or more shorter paths. For example, in the first example above the cross symbol between 2 and 3 in the first Hamiltonian path and cross symbol between 1 and 4 in the second Hamiltonian path, respectively, indicates that the first Hamiltonian path has been broken into two shorter paths: paths 1 2 and 3 4, and the second Hamiltonian path has been broken into shorter paths: paths 31 and 4 2, respectively, forming a path cover:

12
34
31
42

Again, this path cover is not minimal as we can convert this to a minimal path cover constituting 2 paths by combining 1$^{st}$ with last and 2$^{nd}$ with 3$^{rd}$ respectively as follows:

1 2 4
1 3 4

Let us note down some observations related to Hamiltonian and pseudo Hamiltonian path tables, associated with $K_{2n}$ and other graphs containing $2n$ vertices, respectively:

a) Since the degree of each vertex in $K_{2n}$ is $(2n-1)$ therefore every number from 1 to $2n$ that denotes a vertex label must appear $(n-1)$ times as an intermediate vertex label and only one times as an end vertex label in the path table.

b) Every row $j \ \alpha_2^j \ \cdots \ \alpha_j^j \ \cdots \ (2n-j+1)$ of the path table represents a Hamiltonian path incorporating edges $(j, \alpha_2^j), (\alpha_2^j, \alpha_3^j), \cdots$ etc. in succession.

c) Every row $j \ \alpha_2^j \ \cdots \ \alpha_j^j \ \cdots \ (2n-j+1)$ of the path table contains every number from 1 to $2n$ and each number appears only once, since it is a Hamiltonian path.

d) The Hamiltonian path table, $P$, moves to a new Hamiltonian path table under the action of any permutation, $\sigma \in S_{2n}$, where $S_{2n}$ denotes the group of permutations on $2n$ symbols, i.e. under this action $P \rightarrow P^*$ where $P^*$ is also a Hamiltonian path table, containing some different but again Hamiltonian paths, which together still form a Hamiltonian path cover.

e) Every connected graph on at most $2n$ vertices being a subgraph of $K_{2n}$, we can introduce cross symbols "×" in between those numbers representing vertex labels which are not adjacent in the given connected graph, in the Hamiltonian paths where they occur in succession in the path table and construct the pseudo path table.

f) When the broken edges are present only on end sides of a Hamiltonian path in the associated pseudo Hamiltonian path table then the path that is left behind after ignoring the broken edges is still a single path (though not a Hamiltonian path).

g) When broken edges are present only on end sides of each Hamiltonian path in the associated pseudo Hamiltonian path table then all the paths that are left behind after ignoring the broken edges are still a single paths (though not a Hamiltonian paths) and together they form path decomposition made up of paths, less than or equal to $n$ in number, for this graph on $2n$ vertices.

h) For a graph containing 2n vertices if there exists a permutation, $\sigma \in S_{2n}$, such that under its action the associated pseudo Hamiltonian path table, *P*, changes into a new pseudo Hamiltonian path table, *P**, such that all the broken edges now move to end sides, then all the paths that are left behind after ignoring the broken edges are still single paths (though not a Hamiltonian paths) and together they form path decomposition made up of paths, less than or equal to *n* in number, for this graph.

i) For a graph containing 2n or (2n-1) vertices there is a minimal set of permutations formed by selecting them from the set of all possible (2n)! or (2n-1)! permutations representing paths such that the broken edges are present only on end sides of each Hamiltonian path in the associated pseudo Hamiltonian path table, so that when broken edges are removed the left out paths (or truncated permutations) form the minimal path cover for the graph.

j) The pseudo Hamiltonian path cover associated with a graph on 2n vertices tells us that due to cross symbols present between the vertex labels cause breaking of the Hamiltonian paths in the original Hamiltonian path table for $K_{2n}$ have now broken up into paths of shorter lengths or into truncated permutations such that the first Hamiltonian path may have now split into $L_1$ subpaths or truncated permutations, second Hamiltonian path into $L_2$ subpaths or truncated permutations, ….., the last Hamiltonian path into $L_n$ subpaths.

We now proceed to discuss in brief an algorithm to generate disjoint path cover for given connected graph of order 2n or (2n-1).

**Algorithm 5.1:**

1) Construct Hamiltonian path table of order 2n.
2) Introduce cross symbols " × " between the vertex labels which are **not** adjacent in the given connected graph G and thus construct pseudo Hamiltonian path table.
3) Make a list of subpaths formed due to introduction of cross symbols between the vertices appeared in the Hamiltonian paths of the

Hamiltonian path table formed in 1) which are actually not adjacent in the given graph.
4) Use this list to form maximally extended peripheral paths, directly by joining list of subpaths, till every edge of given connected graph belongs to some such path.

**Examples:** 1) We take $G = K_5$. For this graph we consider Hamiltonian path table for $K_6$ as per 1) introducing cross symbols we construct corresponding pseudo Hamiltonian path table as per 2), and we make the list of subpaths as per 3). This list of subpaths is as follows:

        1 2 3 4 5
        3 1 5 2
        5 3
        1 4 2

As per 4) we get the following list of maximally extended peripheral paths which together forms the desired path decomposition for $G = K_5$ as follows:

        1 2 3 4 5
        5 3 1 4 2
        1 5 2

**Problem:** Using observation i) mentioned above show that the minimal path cover thus obtained satisfies Gallai's bound.

## Acknowledgements

I am thankful to Prof. M. R. Modak and Riko Winterle for useful discussions and encouragement.## References

1. L. Lovász, On covering of graphs, in: P. Erdös, G. Katona (Eds.), Theory of Graphs, Academic Press, New York, 1968, pp. 231–236.
2. L. Pyber, Covering the edges of a connected graph by paths, J. Combin. Theory Ser. B 66 (1996) 152–159.
3. Genghua Fan, Path decompositions and Gallai's conjecture, Journal of Combinatorial Theory, Series B 93 (2005) 117 – 125.23